\documentclass[12pt]{article}
\usepackage{amssymb}
\usepackage{amsmath}
\newtheorem{theo}{Theorem}[section]

\newtheorem{prop}[theo]{Proposition}
\newtheorem{lemm}[theo]{Lemma}
\newtheorem{coro}[theo]{Corollary}
\newtheorem{rema}[theo]{Remark}

\def \C{{\mathbb C}}

\def \Z{{\mathbb Z}}
\def \P{{\mathbb P}}
\def \p{{\mathfrak p}}
\def \Q{{\mathbb Q}}
\def \N{{\mathbb N}}

\def \A{{\mathbb A}}
\def \q{{\mathfrak q}}

\begin{document}

\title{Remarks on endomorphisms and rational points}

\author{E. Amerik, F. Bogomolov, M. Rovinsky}

\date{}

\maketitle

\begin{abstract}
Let $X$ be an algebraic variety and let $f: X \dasharrow X$ be a
rational self-map with a fixed point $q$, where everything is defined over
a number field $K$. We make some general remarks concerning the possibility 
of using the behaviour of $f$ near $q$ to produce many rational points on $X$. 
As an application, we give a simplified proof of the potential density of 
rational points on the variety of lines of a cubic fourfold, originally
proved by Claire Voisin and the first author in 2007.
\end{abstract}

(secondary).


\section{Introduction}

Let $X$ be an algebraic variety defined over a number field $K$. One says
that the rational points are potentially dense in $X$, or that $X$ is
potentially dense, if there is a
finite extension $L$ of $K$, such that $X(L)$ is Zariski-dense. For
instance, unirational varieties are obviously potentially dense.
A well-known conjecture of Lang affirms that a variety of general type
cannot be potentially dense; more recently, the question of geometric
characterization of potentially dense varieties has been raised by several
mathematicians, for example, by Abramovich and Colliot-Th\'el\`ene and 
especially by
Campana \cite{C}.
According to their points of view, one expects that the varieties 
with trivial
canonical class should be potentially dense. This is well-known for
abelian varieties, but the simply-connected case remains largely
unsolved. 

Bogomolov and  Tschinkel \cite{BT}  proved potential
density of rational points for $K3$ surfaces admitting an elliptic pencil, 
or an infinite automorphisms group. Hassett and Tschinkel \cite{HT} 
did this for certain symmetric powers of general $K3$-surfaces with a
polarization of a suitable degree. The
key observation of their work was that those symmetric powers are rationally 
fibered in 
abelian varieties over a projective space, and, as the elliptic $K3$ surfaces
of \cite{BT}, they admit a potentially dense multisection which one
can translate by suitable fiberwise rational self-maps 
to obtain the potential density of the ambient
variety. 

More recently, Amerik and Voisin \cite{AV} gave a proof of potential density
for the variety of lines of a sufficiently general cubic fourfold $V$
defined over a number field. Such a variety $X$ is 
an irreducible holomorphic symplectic 
fourfold with cyclic Picard group, in fact the only (up to now) example 
of a simply-connected variety, defined over a number field, 
with trivial canonical class and cyclic Picard
group where the potential density is established. The starting idea
is similar: as noticed in \cite{V1}, $X$ admits a
rational self-map of degree 16. Moreover, $X$ carries a two-parameter
family $\Sigma_b, b\in B$ of surfaces birational to abelian surfaces;
each $\Sigma_b$ parameterizes lines contained in a hyperplane section
of $V$ with three double points.
It is proved in \cite{AV} that under certain genericity conditions
on the pair $(X, b)$, satisfied by many pairs defined over a number field,
the iterates $f^n(\Sigma_b), n\in \N$ are Zariski-dense in $X$. Since
rational points are potentially dense on abelian varieties, this
clearly implies that $X$ is potentially dense. The
proof is rather involved: even the proof of the fact that the number 
of iterates is infinite for some $\Sigma_b$ defined over a number
field is highly non-trivial, using, for instance, $l$-adic Abel-Jacobi
invariants in the continuous \'etale cohomology.

Let us also mention that for varieties of lines of some special cubic
fourfolds (roughly, of those admitting a hyperplane section with six 
double points), the potential density is shown in \cite{HT2} using
the existence of an infinite order birational {\it automorphism} on such
varieties.

The purpose of this article is to further investigate the connection 
between the existence of ``sufficiently nontrivial'' rational self-maps and
the potential density of rational points. In the first part we prove,
among other facts,
that if the differential of a rational self-map $f:X\dasharrow X$ at a 
non-degenerate fixed point $q\in X(\bar \Q)$ has
multiplicatively independent eigenvalues, then the rational points are
potentially dense on $X$. More precisely, we show that under this condition,
one can find a point $x\in X(\bar \Q)$ such that the set of its iterates
is Zariski-dense. Note that this remains unknown for $X$ and $f$ as in 
\cite{AV},
though the question has been raised in \cite{AC} where it is shown
that the map $f:X\dasharrow X$ does not preserve any
rational fibration and therefore the set of the iterates of a general
complex point of $X$ is Zariski-dense.

Unfortunately, it seems to be difficult to find interesting examples
with multiplicatively independent eigenvalues of the differential at a
fixed point. There is certainly plenty of such self-maps on rational
varieties, but since for those the potential density is obvious, we cannot
consider their self-maps as being ``interesting''. In the case of \cite{AV}, 
the eigenvalues of $Df_q$ at a
fixed point $q$ are far from being multiplicatively independent
(lemma \ref{eigenvalues}). We do not know whether the multiplicative 
independence
could hold for a fixed point of a power of $f$.

Nevertheless, even the independence of certain eigenvalues gives interesting
new information. To illustrate this, we exploit this point of view
in the second part, where a simplified proof of the potential
density of the variety of lines of the cubic fourfold is given.

While finishing the writing of this paper, we have learned of a recent 
article \cite{GT},
where a question of a different flavour is approached by similar methods.

\medskip

{\bf Acknowledgements:}
The starting ideas for this article have emerged during the special 
trimester ``Groups and Geometry'' at Centro Di Giorgi, Pisa,  
where the three authors were
present in October 2008. The work has been continued while the first
and the third author were enjoying the hospitality of the 
Max-Planck-Institut f\"ur Mathematik in Bonn, and the second author 
of the IHES, Bures-sur-Yvette. The final version
has been written when the first and the third authors were members 
of the Institute for
Advanced Study in Princeton, supported
by the NSF grant DMS-0635607; the stay of the first author was also 
supported by the Minerva Research Foundation. The second author was 
partially supported by the NSF grant DMS-0701578. 
We would like to thank all these institutions.

We are grateful to S. Cantat, D. Cerveau, A. Chambert-Loir, T.-C. Dinh,
M. Temkin and W. Zudilin for very helpful discussions.

\section{Invariant neighbourhoods}

Let $X$ be a smooth projective variety of dimension $n$ and let 
$f:X\dasharrow X$ be a
rational self-map, both defined over a "sufficiently large" number field $K$.  
We assume that $f$ has a fixed point $q\in X(K)$.
This assumption is not restrictive if, for example, $f$ is a regular polarized
(that is, such that $f^*L=L^{\otimes k}$ for a certain ample 
line bundle 
$L$ and an integer $k>1$) endomorphism: indeed, in this case the set 
of periodic points in
$X(\bar \Q)$ is even Zariski-dense \cite{F}, so replacing $f$
by a power and taking a finite extension of $K$ if necessary, 
we find a fixed point. 

For a number field $K$ we denote by ${\mathcal O}_K$ the ring of integers
of $K$; for a point $\p$, i.e. an equivalence class of valuations of $K$,
$K_{\p}$ denotes the corresponding completion, ${\mathcal O}_{\p}$ the ring 
of integers in $K_{\p}$.

Our starting point is that, for any fixed point $q\in X(K)$ and
a suitable prime ideal $\p\subset {\cal O}_K$, 
we can find a "$\p$-adic neighbourhood" $q\in O_{\p,q}\subset X(K_{\p})$,
on which $f$ is defined and which is $f$-invariant. 

More precisely, choose an affine neighbourhood $U\subset X$ of $q$,
such that the restriction of $f$ to $U$ is regular. By Noether
normalisation lemma, there is a finite $K$-morphism 
$\pi =(x_1,\dots,x_n):U\longrightarrow\A^n_K$ to the affine space, 
which is \'etale at $q$ and which maps $q$ 
to 0. Then the $K$-algebra ${\mathcal O}(U)$ is integral over 
$K[x_1,\dots,x_n]$, 
i.e., it is generated over $K[x_1,\dots,x_n]$ by some regular functions 
$x_{n+1},\dots,x_m$ integral over $K[x_1,\dots,x_n]$. 
The coordinate ring of $U$ is included into the local ring of $q$ and 
the latter is included into its completion: ${\mathcal O}(U)\subset
{\mathcal O}_{U,q}\subset\widehat{{\mathcal O}}_{U,q}=K[[x_1,\dots,x_n]]$. 
In particular, $x_{n+1},\dots,x_m$ become elements of $K[[x_1,\dots,x_n]]$. 
As $f^*$ defines an endomorphism of the ring ${\mathcal O}_{U,q}$ and of 
its completion, the functions $f^{\ast}x_1,\dots,f^{\ast}x_m$ also become 
power series in $x_i$ with coefficients in $K$.

We use the following well-known result (a stronger version
for $n=1$ goes back to Eisenstein) to deduce that the coefficients of 
the power series $x_{n+1},\dots,x_m,f^{\ast}x_1,\dots,f^{\ast}x_m$ are 
$\p$-integral for almost all primes $\p$, i.e., 
$$x_{n+1},\dots,x_m,f^{\ast}x_1,\dots,f^{\ast}x_m\in
{\mathcal O}_K[1/N][[x_1,\dots,x_n]]$$ for some integer $N\ge 1$:

\begin{lemm} \label{fg} Let $k$ be a field of characteristic zero and let
$\phi \in k[[x_1,\dots x_n]]$ be a function algebraic over $k(x_1,\dots,x_n)$. 
Then 
$\phi \in A[[x_1,\dots x_n]]$, where $A$ is a finitely generated
$\Z$-algebra.
\end{lemm}

{\it Proof.} Let $F$ be a minimal polynomial of $\phi$ over $k[x_1,\dots,x_n]$,
so $F(\phi)=0$ and $F'(\phi)\neq 0$. Then $F'(\phi)\in\mathfrak{m}^s
\smallsetminus\mathfrak{m}^{s+1}$ for some $s\ge 0$, where $\mathfrak{m}$
is the maximal ideal in $k[[x_1,\dots,x_n]]$.

Denote by $\phi_d$ the only polynomial of degree $<d$ congruent to $\phi$
modulo $\mathfrak{m}^d$. For a formal series $\Phi$ in $x$ and an integer
$m$ denote by $\Phi_{(m)}$ the homogeneous part of $\Phi$ of degree $m$.
Clearly, $F'(\phi_d)_{(m)}$ is independent of $d$ for $d>m$.

We are going to show by induction on $d$ that the coefficients of the 
homogeneous component of $\phi$ of degree $d$ belong to the $\Z$-subalgebra 
in $k$ generated by coefficients of $F$ (as a polynomial in $n+1$ variables), 
by coefficients of $\phi_{s+1}$ and by the inverse of a certain 
(non-canonical) polynomial $D$ in coefficients of $F$ and 
in coefficients of $\phi_{s+1}$.
To define $D$, choose a discrete valuation $v$ of the field 
$k(x_1,\dots,x_n)$ of rank $n$ trivial on $k$, say such that 
$v(x_i)>v(k(x_1,\dots,x_{i-1})^{\times})$ for all $1\le i\le n$ (equivalently, 
$0<v(x_1^{m_1})<\dots<v(x_n^{m_n})$ for all $m_1,\dots,m_n>0$). 
Then $D$ is the coefficient of the monomial in $F'(\phi_{s+1})_{(s)}$ 
with the minimal valuation. 

For $d\le s$ there is nothing to prove, so let $d>s$. Let $\Delta_d:=
\phi-\phi_d$, so $\Delta_d\in\mathfrak{m}^d$. Then $0=F(\phi_d+\Delta_d)
\equiv F(\phi_d)+F'(\phi_d)\Delta_d\pmod{\Delta_d^2}$, so in particular, 
$F(\phi_d)_{(d+s)}+F'(\phi_d)_{(s)}(\Delta_d)_{(d)}=0$, or equivalently, 
$F(\phi_d)_{(d+s)}+F'(\phi_{s+1})_{(s)}\phi_{(d)}=0$, and thus, 
$\phi_{(d)}=-\frac{F(\phi_d)_{(d+s)}}{F'(\phi_{s+1})_{(s)}}$. 

The field of rational 
functions $k(x_1,\dots,x_n)$ is embedded into its completion with
respect to $v$, and by our choice of $v$ this completion can be identified with
the field of iterated Laurent series 
$k((x_1))\cdots((x_n))$. In particular, $(F'(\phi_{s+1})_{(s)})^{-1}$ 
becomes an iterated Laurent series, whose coefficients are polynomials 
over ${\mathbb Z}$ in the coefficients of $F'(\phi_{s+1})_{(s)}$ and in 
$D^{-1}$ (write $F'(\phi_{s+1})_{(s)}$ as a product of its minimal valuation 
monomial and a rational function, then write the inverse of the
rational function as a geometric series). Then, by induction assumption, the coefficients of $\phi_{(d)}$ 
are polynomials over ${\mathbb Z}$ in coefficients of $F$, 
coefficients of $\phi_{s+1}$ and in $D^{-1}$.

\

Therefore, for almost all primes $\mathfrak p\subset {\cal O}_K$,
the coefficients of the power series 
$x_{n+1}, \dots, x_{m},f^*x_1, \dots f^*x_n$ are 
integral in $K_{\p}$. Choose a $\p$ satisfying this condition and such
that, moreover,  the irreducible polynomials $P_i(X)=
P_i(x_1,\dots,x_n;X)\in K[x_1,\dots,x_n;X]$ which are minimal monic 
polynomials of $x_i$ for $n<i\leq m$, have $\p$-integral coefficients and
the elements $P_i'(0,\dots,0;x_i(q))$ are invertible 
in ${\mathcal O}_{\p}$ for each $n<i\leq m$ (this last condition holds
for almost all primes $\p$ since the morphism $x$ is \'etale at $q$).

Define the system of $\p$-adic neighbourhoods $O_{{\p}, q, s}$, $s\geq 1$ 
of the point $q$ as follows:
$$O_{{\p}, q, s}=\{t\in U(K_{\p})|x_i(t)\equiv x_i(q)\pmod {\p^s}\ for
\ 1\leq i\leq m\}.$$ We set $O_{\p,q}:=O_{\p,q,1}$.

\begin{prop} \label{nbhd}  (1) The functions $x_1, \dots x_n$ give a bijection
between $O_{{\p}, q, s}$ and the n-th cartesian power of $\p^s$.

(2) The set $O_{\p,q}$ contains no indeterminacy points of $f$.

(3) $f(O_{{\p}, q, s})\subset O_{{\p}, q, s}$ for $s\geq 1$. Moreover, 
$f:O_{\p,q,s}\stackrel{\sim}{\longrightarrow}O_{\p,q,s}$ is bijective 
if $\det Df_q$ is invertible in ${\mathcal O}_{\p}$.

(4) The $\bar \Q$-points are dense in $O_{{\p}, q, s}$.

\end{prop}

{\it Proof:} These properties are clear from the definition and 
the inclusion of the elements $x_1,\dots,x_m,f^{\ast}x_1,\dots,
f^{\ast}x_m$ into ${\mathcal O}_{\p}[[x_1,\dots,x_n]]$. 
\begin{enumerate} \item The map $x$ from $O_{\p,q,s}$ to the $n$-th 
cartesian power of $\p^s$ is injective, since the coordinates 
$x_{n+1},\dots,x_m$ of a point $t$ are determined uniquely by the 
coordinates $x_1,\dots,x_n$ and the condition $t\in O_{\p,q}$. 
Let $P_i(x_i)=P_i(x_1,\dots,x_n;x_i)=0$ be the minimal monic polynomial 
of $x_i$ for $n<i\le m$. For fixed values of $x_1,\dots,x_n\in\p$ this 
equation has precisely $\deg P_i$ solutions (with multiplicities) in 
$\bar K_{\p}$. As $P_i'(x_i(q))\in{\mathcal O}_{\p}^{\times}$, 
$x_i(t)\equiv x_i(q)\pmod{\p}$ is a simple root of $P_i\pmod{\p}$, 
and thus, any root congruent to $x_i(q)$ modulo $\p$ 
is not congruent to any other root modulo $\p$. 

It is surjective, since $x_{n+1},\dots,x_m$ are 
convergent series on $\p^s$ with constant values modulo $\p^s$. 
\item The functions $f^{\ast}x_i$, $1\le i\le m$, are convergent series 
on $O_{\p,q}$. 
\item The functions $f^{\ast}x_i$ on the $n$-th cartesian power of $\p^s$, 
$1\le i\le m$, are constant modulo $\p^s$. This shows that 
$f(O_{\p,q,s})\subseteq O_{\p,q,s}$. If $\det Df_q\neq 0$ then the inverse 
map $f^{-1}$ is well-defined in a neighbourhood of 0. If 
$\det Df_q\in{\mathcal O}_{\p}^{\times}$ then $f^{-1}$ is defined by series in 
${\mathcal O}_{\p}[[x_1,\dots,x_n]]$, i.e., it is well-defined on $O_{\p,q}$. 
\item The $\bar\Q$-points are dense in the $n$-th cartesian power of 
$\p^s$ and they lift uniquely to $O_{\p,q,s}$.
\end{enumerate}

\medskip

Let $\lambda_1, \dots, \lambda_n$ be the eigenvalues of the tangent map
$Df_q$. We assume that $q$ is a non-degenerate fixed point of $f$, so that
$\lambda_i\neq 0$. Note that the $\lambda_i$ are algebraic numbers.
Extending, if necessary, the field $K$, we may assume
that $\lambda_i\in K$. The following is a consequence of the $p$-adic 
versions of 
several well-known results in
dynamics and number theory:

\begin{prop}\label{linear-gen}  Assume that $\lambda_1, \dots, \lambda_n$ are multiplicatively
independent. Then in some $\p$-adic neighbourhood $O_{{\p}, q, s}$, 
the map $f$ is
equivalent to its linear part $\Lambda$ (i.e. there exists a formally 
invertible $n$-tuple of formal power series $h=(h^{(1)},\dots,h^{(n)})$ 
in $n$ variables $(x_1,\dots,x_n)=x$, convergent together 
with its formal inverse on a neighbourhood of zero, such that 
$h(\lambda_1x_1,\dots,\lambda_nx_n)=f(h(x_1,\dots,x_n))$.

\end{prop}

{\it Proof}: It is well-known that in absence of relations 
$$\lambda_1^{m_1}...\lambda_n^{m_n}=\lambda_j,\ 1\leq j\leq n,\ m=\sum m_i\geq 2,
\ m_i\geq 0$$ ("resonances"), there is a unique
formal linearization of $f$, obtained by formally solving the 
equation $f(h(x))=h(\Lambda(x))$; the expressions $\lambda_1^{m_1}...\lambda_n^{m_n}-\lambda_j$
appear in the denominators of the coefficients of $h$ (see for example \cite{Arn}). 
The problem
is of course whether $h$ has non-zero radius of convergence, that is, 
whether the denominators are ``not too small''.
By Siegel's theorem (see \cite{HY} for its $p$-adic version) this is
the case as soon as the numbers $\lambda_i$ satisfy the diophantine
condition $$|\lambda_1^{m_1}...\lambda_n^{m_n}-\lambda_j|_p>Cm^{-\alpha}$$
for some $C, \alpha$. By \cite{Yu} , this condition is always satisfied by 
algebraic numbers.

\medskip

When the fixed point $q$ is not isolated, the eigenvalues $\lambda_i$
are always resonant. However, as follows from the results proved in the
Appendix, if ``all resonances come from the fixed subvariety'', the 
linearization is still possible.

More precisely, extending $K$ if necessary, we may choose  
$q\in X(K)$ which is a smooth point of the fixed point locus 
of $f$. Let $F$ be the
irreducible component of this locus containing $q$. Let $r=\dim F$.
By a version of Noether's normalization lemma 
\cite[Theorem 13.3 \& geometric interpretation on p.284]{eisen}, 
we may assume that our finite morphism $\pi=(x_1,\dots,x_n):U\longrightarrow
\A^n_K$ which is \'etale at $q$, maps $F\cap U$ onto the coordinate plane
 $\{x_{r+1}=\dots=x_n=0\}$.

\begin{prop}\label{linear-still}  
Let $q$ be a general fixed point of $f$ as above. 
Suppose that the tangent map $Df_q$ is semisimple and that its 
eigenvalues $\lambda_1,\dots,\lambda_n$ satisfy the condition 
$\lambda_{r+1}^{m_{r+1}}\cdots\lambda_n^{m_n}\neq\lambda_i$ for 
all integer $m_{r+1},\dots,m_n\ge 0$ with $m_{r+1}+\dots+m_n\ge 2$ 
and all $i$, $r<i\le n$ (and $\lambda_1=\dots=\lambda_r=1$). 

Suppose that the eigenvalues of $Df_q$ do not vary with 
$q$. Then, for each $\p$ as above, the map $f$ can be linearized in some 
$\p$-adic neighbourhood $O_{\p,q,s}$ of $q$, i.e., there exists a formally 
invertible $n$-tuple of formal power series $h=(h^{(1)},\dots,h^{(n)})$ 
in $n$ variables $(x_1,\dots,x_n)=x$ convergent together 
with its formal inverse on a neighbourhood of zero such that 
$h(\lambda_1x_1,\dots,\lambda_nx_n)=f(h(x_1,\dots,x_n))$. \end{prop}

{\it Proof:} Taking into account that the $\lambda_i$ are algebraic
numbers and so the diophantine condition 
$$|\lambda_{r+1}^{m_{r+1}}\cdots\lambda_n^{m_n}-\lambda_j|_p
>C(\sum_{r<i\le n}m_i)^{-\alpha}$$ 
is automatically satisfied for some $C,\alpha>0$, whenever 
$\sum_{r<i\le n}m_i\ge 2$, this is just Theorem \ref{linearizaciq}. 

\begin{rema}\label{equalok} It might be worthwhile to mention explicitely 
that since the sum of the $m_i$ in the condition of this proposition
must be at least two, the relation $\lambda_i=\lambda_j$ is NOT a resonance 
for 
$i,j$ between $r+1$ and $n$, in other words, linerization is still
possible when some (or even all) of the non-trivial eigenvalues are equal.
\end{rema}

In order to apply these propositions to the study of iterated orbits
of algebraic points, we need the following lemma.

\begin{lemm} \label{topol-lin-nezav-harakt} Let $a_1,a_2,\dots$ be a 
sequence in $K_{\p}^{\times}$ tending to 0. Let $b_1,b_2,\dots$ be a sequence 
of pairwise distinct elements in 
${\cal O}_{\p}^{\times}$ generating a torsion-free subgroup. 
Then any infinite subset $S\subset{\mathbb N}$ 
contains 
an element $s$ such that $\sum_{i\ge 1}a_ib_i^s\neq 0$ (in other words, the
set of $s\in \N$ such that $\sum_{i\ge 1}a_ib_i^s = 0$ is finite). 
\end{lemm}

{\it Proof:} Renumbering if necessary, we may suppose that
$|a_1|=|a_2|=\dots =|a_N|>|a_i|$ for any $i>N$. Suppose that 
$\sum_{i\ge 1}a_ib_i^s = 0$ for every $s\in S$. 

First assume that $S={\mathbb N}$. It follows that for any polynomial
$P$, $\sum_{i\ge 1}a_iP(b_i) = 0$. By the triangular inequality, we'll
get a contradiction as soon as we find a polynomial $P$ such that
$|P(b_k)|<|P(b_1)|$ for $2\leq k\leq N$ and $|P(b_k)|\leq|P(b_1)|$ for $k>N$.
To construct such a $P$, choose an ideal ${\q}={\p}^s$ such that the 
$b_i$ are pairwise distinct modulo ${\q}$ for $i=1,\dots, N$ and let 
$$P(x)=\prod_{b\in {\cal O}_{\p}/{\q},\ b_1\not\in b}(x-\bar b),$$
where $\bar b$ denotes any representative of the class $b$.
An easy check gives that $|P(x)|=|P(b_1)|$ when $x\equiv b_1 \pmod \q$ and 
$|P(x)|<|P(b_1)|$ otherwise, so $P$ has the required properties.

(This polynomial $P$ has been indicated to us by A. Chambert-Loir.)

Now let $S\subset \N$ be arbitrary. Take an integer $M$ such that
$|b_i^M-1|<|p|^{\frac{1}{p-1}}$ for all $i$ (this is possible since $b_i$
are in ${\cal O}_{\p}^{\times}$ and the residue field is finite), then
$c_i=\log b_i^M$ is defined for all $i$ and $\exp (c_i)=b_i^M$. 
Since the subgroup generated by the $b_i$ is torsion-free, all $c_i$ are 
different. We claim that for certain coefficients $a'_i$, we can
write identities $\sum_{i\ge 1}a'_ic_i^m = 0$ for all $m\in \N$, so 
this case reduces to that of $S=\N$.

Indeed, consider $S$ as a subset of $\Z_p$; it has a limit point $s_0$.
For a sequence of $m$-tuples $j_1>j_2>\dots>j_m$
in $S$ such that $j_1\equiv j_2\equiv\dots\equiv j_m\pmod M$ and 
$\lim j_1=\lim j_2=\dots=\lim j_m=s$, and any analytic function $g$
on $\Z_p$, one has
$$\frac{g^{(m)}(s_0)}{m!}=
\lim\sum_{l=1}^m\frac{g(j_l)}{\prod_{k\neq l}(j_l-j_k)}$$
(Newton interpolation formula).

Pick a class modulo $M$ which contains a sequence in $S$ converging to
$s_0$. There is an analytic function $g_i$ such that $g_i(k)=b_i^k$
when $k$ is in this class. By the formula above, we have
$$g_i^{(m)}(s_0)=\frac{c_i^m}{M^m}g_i(s_0)=
m!\lim\sum_{l=1}^m\frac{g(j_l)}{\prod_{k\neq l}(j_l-j_k)},$$
where all the $j_l$ are in the same class modulo $M$.
This gives $$\sum_{i=1}^{\infty}a_ig_i(s_0)c_i^m=
M^m m!\lim \sum_{l=0}^m\frac{1}{\prod_{k\neq l}(j_l-j_k)}\sum_{i=1}^{\infty}
a_ig_i(j_l)=0$$ since the $j_l$ are in $S$, q.e.d.

\medskip

From now on, we assume that $\p$ is chosen such that all $\lambda_i$
belong to ${\cal O}^{\times}_{\p}$ (this is of course the case for
almost all $\p$.

The first part of the lemma (case $S=\N$) immediately implies
the following corollary:

\begin{coro} \label{indep-dense} If $\lambda_1, \dots, \lambda_n$ are multiplicatively independent,
the rational points on $X$ are potentially dense. More precisely, there is
a point $x\in X(\bar \Q)$ with Zariski-dense iterated orbit 
$\{f^i(x)|i\in \N\}$.
\end{coro}

{\it Proof:} Since algebraic points are dense in $O_{{\p}, q, s}$, we can find a point
$x\in X(\bar \Q)$ which is contained in $O_{{\p}, q, s}$, away from the coordinate
hyperplanes in the local coordinates $(y_1,\dots y_n)$ linearizing $f$. We claim that the
iterated orbit of such a point is Zariski dense in $X$. Indeed, if not,
there is a regular function $G$ on $U$ vanishing on $f^i(x)$ for all $i$; in the
local linearizing coordinates on $O_{{\p}, q, s}$, $G$ becomes a
convergent power series $G=\sum_I a_Iy^I$. If $x=(x_1, \dots x_n)$,
we get $\sum_Ia_Ix^I(\lambda^I)^i=0$ for all $i\in \N$. Since the $\lambda_i$
are multiplicatively independent, the numbers $\lambda^I$ are distinct,
contradicting lemma \ref{topol-lin-nezav-harakt}.

\medskip
 
Another useful version of this corollary is the following:
\begin{coro} \label{analytic-closure} Let $n$ be a natural number, 
${\cal T}$ be the $n$-th cartesian power of ${\mathcal O}^{\times}_{\p}$ and 
$\Lambda=(\lambda_1,\dots,\lambda_n)\in {\cal T}$, where $\lambda_1,\dots,\lambda_n$ 
are multiplicatively independent. Let $S\subset{\mathbb N}$ be an 
infinite subset. Then the set $\{\Lambda^i~|~i\in S\}$ is ``analytically 
dense'' in ${\cal T}$, i.e., for any non-zero Laurent series $F=\sum_Ia_Iy^I$ 
convergent on ${\cal T}$, there is an $s\in S$ such that $F(\Lambda^s)\neq 0$. 
\end{coro}
{\it Proof.} Otherwise, after a renumbering of $a_I$ as $a_i$ and setting 
$b_i=\lambda^I$ for $i$ corresponding to $I$, we get 
$\sum_{i\ge 1}a_ib_i^s=0$ for all $s\in S$. This contradicts Lemma 
\ref{topol-lin-nezav-harakt}. 

\medskip

The case $S\neq \N$ yields the following remark which might be useful
in the study of the case when the fixed point $q$ is not isolated
(as in the next section).

\begin{coro} \label{irreducibility} Under assumptions of Proposition 
\ref{linear-still}, let $Y\subset U$ be such an irreducible subvariety that,
possibly after a finite field extension,  
$Y(K_{\p})$ meets a sufficiently small $\p$-adic neighbourhood of $q$. 
Suppose that the multiplicative subgroup $H\subset{\mathcal O}_{\p}^{\times}$ 
generated by $\lambda_1,\dots,\lambda_n$ is torsion-free. Let 
$S\subset{\mathbb N}$ be an infinite subset. Then the Zariski closure of 
the union $\bigcup_{i\in S}f^{i}(Y)$ is independent of $S$, and 
therefore, is irreducible. \end{coro}
{\it Proof.} To check the independence of $S$, let us show that any regular 
function $F$ on $U$ vanishing on $\bigcup_{i\in S}f^{i}(Y)$ vanishes 
also on $\bigcup_{i\ge 1}f^{i}(Y)$. 

Let $y_1,\dots,y_n$ be local coordinates at $q$ linearizing and 
diagonalizing $f$ on a neighbourhood $O_{\p,q,s}$. Then for any choice of
a sufficiently large number field $K$, 
$f^{\circ i}(Y)(K_{\p})\cap O_{\p,q,s}$ contains 
$f^{\circ i}(Y(K_{\p})\cap O_{\p,q,s})$, and therefore, each irreducible 
component of the Zariski closure of $\bigcup_{i\in S}f^{i}(Y)$ 
meets $O_{\p,q,s}$. 

This implies that we can work in $O_{\p,q,s}$, where $F$ becomes 
a convergent power series $\sum_Ia_Iy^I$. For any point 
$t=(t_1,\dots,t_n)\in Y( K_{\p})\cap O_{\p,q,s}$ the series $F$ 
vanishes at the points $f^{i}(t)$ for all $i\in S$ (i.e., 
$\sum_Ia_I\lambda^{iI}t^I=0$ for all $i\in S$) and it remains to show 
that $F$ vanishes also at the points $f^i(t)$ {\it for all} $i\ge 1$. 

Let us work only with multi-indices $I$ such that $t^I$ is non-zero.
The set of such multi-indices splits into the following 
equivalence classes: $I\sim I'$ if $\lambda^I=\lambda^{I'}$. By Lemma 
\ref{topol-lin-nezav-harakt}, the sum of $a_It^I$ over each equivalence class 
is zero, and thus, $\sum_Ia_I\lambda^{iI}t^I=0$ for all $i\ge 1$. 

\medskip

Finally, the following is an obvious generalization of \ref{indep-dense}:

\begin{coro} \label{density-of-orbits} Under the assumptions of Proposition 
\ref{linear-still}, let $t$ be a sufficiently general algebraic point of $U$ 
in a sufficiently small $\p$-adic neighbourhood of $q$ and $H$ be the 
multiplicative group generated by $\lambda_1,\dots,\lambda_n$. Then the 
dimension of the Zariski closure of the $f$-orbit $\{f^{i}(t)~|~
i\in{\mathbb N}\}$ is greater than or equal to $r=rank(H)$. 
\end{coro}
{\it Proof.} Let $y_1,\dots,y_n$ be local coordinates at $q$ linearizing and 
diagonalizing $f$ in a neighbourhood $O_{\p,q,s}$. Replacing $f$ by a 
power, we may assume that the eigenvalues $\lambda_i$ generate
a torsion-free group. Take a point $t=(t_1,\dots, t_n)\in O_{\p,q,s}$ away from
the coordinate hyperplanes and consider the natural embedding
of the $n$-th cartesian power ${\cal T}$ of ${\cal O}^{\times}_{\p}$:
$$t: {\cal T}\to O_{\p,q,s},\ (\mu_1,\dots,\mu_n)\mapsto 
(\mu_1t_1,\dots,\mu_n t_n).$$
On the image of ${\cal T}$, we can find a monomial coordinate change
such that in the new coordinates $y'_1,\dots,y'_n$ one has
$$f(y'_1,\dots, y'_n)=
(\lambda'_1y'_1,\dots, \lambda'_ry'_r, y'_{r+1},\dots y'_n)$$
and $\lambda'_i$, $1\leq i\leq r$ are multiplicatively independent
(Euclid's algorithm). Now apply Corollary \ref{analytic-closure}
(notice that since it is about Laurent series, we can allow monomial
coordinate changes as above).

\medskip

To end this section, let us mention that our considerations yield a
simple and self-contained proof of the well-known result that rational
points are potentially dense on abelian varieties (see for example
\cite{HT}, Proposition 3.1). The following is a result of a discussion
with A. Chambert-Loir:
\begin{prop}\label{abvar} Let $A=A_1\times \dots \times A_n$ be a product
of simple abelian varieties, let $k_1, \dots k_n$ be multiplicatively
independent positive integers and set $f: A\to A, f=f_1\times\dots \times f_n,$
where for each $i$, $f_i$ is the multiplication by $k_i$ on $A_i$.
Then there exist points $x\in A(\bar \Q)$ with Zariski-dense iterated orbit
$\{f^m(x)|m\in \N\}$.
\end{prop}
{\it Proof:} Take a suitable number field $K$, $\p\subset {\cal O}_K$ and 
for each $i$, a $\p$-adic neighbourhood $O_i$ such that $f_i$ is 
linearized in $O_i$. Work in the linearizing coordinates. The same argument
as in the above corollaries shows that for a
 point $(x_1, \dots , x_n)\in O_1\times \dots \times O_n$, $x_i\neq 0$,
the iterated orbit is analytically dense in $l_1\times \dots \times l_n$,
where $l_i\subset O_i$ is the line through the origin generated by $x_i$.
Induction by $n$ and the simplicity of the $A_i$ then show that 
$l_1\times \dots \times l_n$ is not contained in a proper abelian subvariety
of $A$. Since its Zariski closure is invariant by $f$, it must coincide
with $A$.

\medskip

Since any abelian variety is isogeneous to a product of simple abelian 
varieties, this proposition proves that abelian varieties are potentially
dense. Note that, unlike Proposition 3.1 from \cite{HT}, it does not
show the existence of an algebraic point generating a Zariski-dense
subgroup.

\begin{rema}\label{zhang} Let $f$ be a polarized endomorphism of a smooth
projective variety $X$. S.-W. Zhang conjectures in \cite{Z} that one can
find a point in $X(\bar \Q)$ with Zariski-dense iterated orbit. According 
to our corollary \ref{indep-dense}, this conjecture is true provided that
some power of $f$ has a non-degenerate fixed point with multiplicatively 
independent 
eigenvalues of the tangent map. Unfortunately not all polarized
endomorphisms have this property: a simple example is the endomorphism
of $\P^n$, $n>1$, taking all homogeneous coordinates to the $m$-th power, $m>1$. 
\end{rema}

\section{Variety of lines of the cubic fourfold}

The difficulty in using the results of the previous section to prove
potential density of rational points is that it can be hard to find 
an interesting
example such that the eigenvalues of the tangent map at some
fixed point are multiplicatively independent. For instance if 
$f$ is an automorphism and $X$ is a projective $K3$ surface, or, more
 generally,
an irreducible holomorphic symplectic variety,
then the product of the eigenvalues is always a root of unity, as noticed
for instance in \cite{Bv}. 

So even when a linearization in the neighbourhood of a fixed point is possible,
the orbit of a general algebraic point may be contained in a relatively small 
analytic subvariety of the neighbourhood (of course this subvariety does not 
have
to be algebraic, but it is unclear how to prove that it actually is not). 
Nevertheless, with some
additional geometric information, one can still follow this approach to
prove the potential density.

In the rest of this note, we illustrate this by giving a simplified
proof of the potential density of the variety of lines of a cubic fourfold, 
which is the main result of \cite{AV}. The proof uses several ideas
from \cite{AV}, but we think that certain aspects become more
transparent thanks to the introduction of the ``dynamical'' point of view and
the use of $\p$-adic neighbourhoods. 

We recall the setting of \cite{AV} (the facts listed below are taken from 
\cite{V1} and \cite{A}).
Let $V$ be a general smooth cubic in $\P^5$ and let
$X\subset G(1,5)$ be the variety of lines on $V$. This is an irreducible
holomorphic symplectic fourfold: $H^{2,0}(X)$ is generated by a nowhere vanishing
form $\sigma$. 
For $l\subset V$ general, there is a unique plane $P$ tangent to $X$
along $l$ (consider the Gauss map, it sends $l$ to a conic in the dual
projective space). The map $f$ maps $l$ to the residual line $l'$. It 
multiplies the form $\sigma$ by $-2$; in particular, its degree is 16. 
The indeterminacy locus $S$ consists of points such that the image of the 
corresponding
line by the Gauss map is a line (and the mapping is 2:1). This is a smooth 
surface of general type, resolved by a single blow-up. For a general $X$,
the Picard group is cyclic and thus the Hodge structure on $H^2(X)^{prim}$
is irreducible (thanks to $h^{2,0}(X)=1$); the space of algebraic cycles
is generated by $H^2=c_1^2(U^*)$ and $\Delta=c_2(U^*)$, where $U$ is the restriction
of $U_{G(1,5)}$, the
universal rank-two bundle on $G(1,5)$. One deduces from Terasoma's 
theorem \cite{T} that these
conditions are satisfied by a ``sufficiently general'' $X$ defined over a number field, in fact
even over $\Q$; ``sufficiently general'' meaning ``outside of a thin subset in
the parameter space''.
One computes that the cohomology class of $S$ is $5(H^2-\Delta)$ to conclude 
that
$S$ is irreducible and non-isotropic with respect to $\sigma$.

\subsection{Fixed points and linearization}

The fixed point set $F$ of our rational self-map $f:X\dasharrow X$ is the set of 
points such that along the corresponding
line $l$, there is a tritangent plane to $V$. Strictly speaking, this is
the closure of the fixed point set, since some of such points are in
the indeterminacy locus; but for simplicity we shall use the term 
``fixed point set'' as far as there is no danger of confusion.

\begin{prop}\label{fixed} The fixed point set $F$ of $f$ is an isotropic surface
of general type.
\end{prop}

{\it Proof:} It is clear from $f^*\sigma=-2\sigma$ that $F$ is isotropic.
Let $I\subset G(1,5)\times G(2,5)$ with projections
$p_1, p_2$ be the incidence
variety $\{(l,P)| l\subset P\}$ and let ${\cal F}\subset I\times 
\P H^0({\cal O}_{\P^5}(3))$ denote the 
variety of triples $\{(l,P,V)| V\cup P=3l\}$. This is a projective bundle
over $I$, so ${\cal F}$ is smooth and thus its fiber $F'_V$ over a general
$V\in \P H^0({\cal O}_{\P^5}(3))$ is also smooth. This fiber clearly
projects generically one-to-one on the corresponding $F=F_V$, since
along a general line $l\subset V$ there is only one tangent plane, 
and a fortiori only one
tritangent plane if any; so $F'=F'_V$ is a desingularization of $F$. 
Since $dim(I)=11$ and since intersecting the plane $P$ along the triple
line $l$ imposes 9 conditions on a cubic $V$, we conclude that 
$F'$ and $F$ are surfaces.

To compute the canonical class, remark that $F'$ is the zero locus of a
section of a globally generated vector bundle on $I$. This vector bundle
is the quotient of $p_2^*S^3U_{G(2,5)}^*$ (where $U_{G(2,5)}$ denotes the tautological 
subbundle on $G(2,5)$) by a line subbundle ${\cal L}_3$ whose fiber at $(l,P)$ is the
space of degree 3 homogeneous polynomials on $P$ with zero locus $l$.
One computes that the class of ${\cal L}_3$ is three times the difference
of the inverse images of the Pl\"ucker hyperplane classes on $G(2,5)$ and $G(1,5)$, and it follows that 
the canonical class of $F$ is $p_2^*(3c_1(U^*))$, which is ample 
(we omit the details since an analogous
computation is given in \cite{V}, and a more detailed version of it 
in \cite{Pc}).

\begin{rema}\label{fixed-indet} Since $F$ is isotropic and $S$ is not,
$S$ cannot coincide with a component of $F$. In fact, dimension count
shows that $F\cap S$ is a curve.
\end{rema}

\begin{prop}\label{eigenvalues} For a general (that is, non-singular and out
of the indeterminacy locus) point $q\in F$, the tangent 
map
$Df_q$ is diagonalized with eigenvalues $1, 1, -2, -2$.
\end{prop}

{\it Proof:} This follows from the fact that $f^*\sigma = -2\sigma$.
and the fact that the map is the identity on the lagrangian plane 
$T_pF\subset T_pX$. Let
$e_1,e_2,e_3,e_4$ be the Jordan basis with $e_1, e_2\in T_pF$.
There is no Jordan cell corresponding to the eigenvalue 1, since
in this case $e_4$ would be an eigenvector with eigenvalue 4,
but then $\sigma(e_1,e_4)=\sigma(e_2,e_4)=\sigma(e_3,e_4)=0$, contradicting 
the fact
that $\sigma$ is non-degenerate.
By the same reason, the eigenvalues at $e_3$ and $e_4$ are both equal to
$\pm 2$. 
Suppose that $Df_q$ is not diagonalized, so sends $e_3$ to $\pm 2e_3$ and 
$e_4$ to $e_3\pm 2e_4$. In both cases $\sigma(e_3, e_4)=0$. If $e_3$ goes to
$2e_3$, we immediately see that $e_3\in Ker(\sigma)$, a contradiction.
Finally, if $Df_q(e_3)=-2e_3$ and $Df_q(e_4)=e_3-2e_4$, 
we have 
$$-2\sigma(e_1,e_4)=\sigma(e_1,e_3)-2\sigma(e_1,e_4),$$
so that $\sigma(e_1,e_3)=0$, but by the same reason $\sigma(e_2,e_3)=0$,
again a contradiction to non-degeneracy of $\sigma$.

\begin{prop}\label{linear} 
(1) Let $q\in X(K)$ be a general fixed point of $f$ as above 
and let $O_{\p,q}$ be its
$\p$-adic neighbourhood for a suitable $\p$, 
as in the previous section. Then $f$ is equivalent to its linear
part in a sufficiently small subneighbourhood $O_{\p,q,s}$; that is,
there exist power series $h=h_q$ in four 
variables $(t_1,t_2,t_3,t_4)=t$ such that $h(t_1,t_2,-2t_3,-2t_4)=f\circ h(t)$,
convergent together with its inverse in some neighbourhood of zero.  

(2) In the complex setting, the analogous statements are true. Moreover,
the maps $h_{t_1,t_2},$ where $h_{t_1,t_2}(x,y)=h(t_1,t_2,x,y)$
extend to global meromorphic maps from $\C^2$ to $X$.

\end{prop}

{\it Proof:} 1) Since the couple of non-trivial eigenvalues $(-2,-2)$
is non-resonant, this is just the Proposition \ref{linear-still}.

2) In the complex case, the linearization is a variant of a classical 
result due to Poincare \cite{P}. One writes the formal power series
in the same way as in the $p$-adic setting, thanks to the absence of
the resonances; but it is much easier to prove its convergence thanks
to the fact that now $\lambda_3=\lambda_4=-2$ and $|-2|>1$, and so
the absolute values of the denominators which appear when one computes
the formal power series are bounded from
below (these denominators are in fact products of the factors of the form
$\lambda_3^{m_3}\lambda_4^{m_4}-\lambda_i$ for $m_3+m_4\geq 2$, 
$m_3, m_4\geq 0$). 
For the sake of brevity,
we refer to \cite{Rg} which proves the analogue of Theorem
\ref{linearizaciq}, case (2), in the complex case and under a weaker
diophantine condition on the eigenvalues (Rong assumes moreover that
$|\lambda_i|=1$, but as we have just indicated, in our case all estimates
only become easier, going back to \cite{P}).

To extend the maps $h_{t_1,t_2}$ to $\C^2$, set $$h_{t_1,t_2}(x)=
f^k(h_{t_1,t_2}((-2)^rx)),$$
where $(-2)^rx$ is sufficiently close to zero; one checks that this is 
independent of choices.

\medskip

We immediately get the following corollary (which follows from the results
of \cite{AV}, but for which there was as yet no elementary proof):

\begin{coro}\label{nonpreppoints}
There exist points in $X(\bar \Q)$ which are not preperiodic for $f$.
\end{coro} 

{\it Proof:} Indeed, $\bar \Q$-points are dense in $O_{\p, q}$. Take one in 
a suitable 
invariant subneighbourhood and use the linearization given by the proposition
above.

\medskip

\begin{rema}\label{heights}
If f were regular, this would follow from the theory of canonical heights;
but this theory does not seem to work sufficiently well for polarized 
rational self-maps.
\end{rema}

\subsection{Non-preperiodicity of certain surfaces}

The starting point of \cite{AV} was the observation that $X$ is covered
by a two-parameter family $\Sigma_b, b\in B$  of birationally abelian 
surfaces, namely,
surfaces parametrizing lines contained in a hyperplane section of $V$
with 3 double points. On a general $X$, a general such surface has cyclic 
Neron-Severi group (\cite{AV}); moreover, many of those surfaces $\Sigma$
defined over a number field have the same property, as shown by an argument
similar to that of Terasoma \cite{T}. 

In \cite{AV}, it is shown that the iterates of a suitable $\Sigma$
defined over a number field and with cyclic Neron-Severi group are 
Zariski-dense in $X$. The first step is to prove non-preperiodicity, 
that is,
the fact that the number of $f^k(\Sigma)$, $k\in \N$, is infinite.
Already at this stage the proof is highly non-trivial, using the $l$-adic
Abel-Jacobi invariant in the continuous \'etale cohomology.

In this subsection, we give an elementary proof of the non-preperiodicity
of a suitable $\Sigma$,
which is based on proposition \ref{linear}. Moreover, this works without an assumption
on its N\'eron-Severi group, and also for any $X$, not only for
a ``general'' one. 

\begin{lemm}\label{noninv} The surface $\Sigma$ is never invariant by $f$.
\end{lemm}

{\it Proof:} The surface $\Sigma$ is the variety of lines contained
in the intersection $Y=V\cap H$, where $H$ is a hyperplane in $\P^5$
tangent to $V$ at exactly three points. For a general line $l$
corresponding to a point of $\Sigma$, there is a unique plane $P$
tangent to $V$ along $l$, and the map $f$ sends $l$ to the residual
line $l'$. If $\Sigma$ is invariant, $l'$ and therefore $P$ lie in $H$,
and $P$ is tangent to $Y$ along $l$. But this means that $l$ is
"of the second type" on $Y$ in the sense of Clemens-Griffiths
(i.e. the Gauss map of $Y\subset H=\P^4$ maps $l$ to a line in
$(\P^4)^*$ as a double covering), see \cite{CG}. At the same time it 
follows from the results
of \cite{CG} that a general line on a cubic threefold with double
points is "of the first type" (mapped bijectively onto a conic by the
Gauss map), a contradiction.

\medskip

Passing to the $p$-adic setting and taking a $\Sigma$ meeting a
small neighbourhood of a general fixed point $q$ of $f$, 
we see by corollary \ref{irreducibility}
that the Zariski closure of $\cup_k f^k(\Sigma)$ is irreducible. 
Since $\Sigma$ cannot be $f$-invariant by the lemma above,
this means that $\Sigma$ is not preperiodic and so the Zariski closure $D$
of $\cup_k f^k(\Sigma)$ is at least a divisor.

\medskip

Coming back to the complex setting and taking a $\Sigma$ passing close
to $q$ in both $p$-adic and complex topologies, let us make a 
few remarks on the
geometry of $D$.

In a neighbourhood of our fixed point $q$, the intersections of
$D$ with the images of $h_{t_1,t_2}$ are $f$-invariant analytic
subsets. From the
structure of $f$ as in Proposition \ref{linear} we deduce that 
such a subset is either the whole image of $h_{t_1,t_2}$, or
a finite union of "lines through the origin" (that is, images
of such lines by $h_{t_1,t_2}$). If the last case holds generically, 
$D$ must
contain $F$ by dimension reasons. If the generic case is the first one,
$D$ might only have a curve in common with $F$.

To sum up, we have the following

\begin{theo}\label{nonprep} For $\Sigma$ meeting a sufficiently small
neighbourhood of $q$, the Zariski closure $D$ of $\cup_kf^k(\Sigma)$ 
is of dimension at least 
three. If it is of dimension three, this is an irreducible divisor which 
either contains the surface of fixed points $F$, or has a curve in common
with $F$. In this last case, $D$ contains correspondent "leaves" 
(images of $\C^2$ from \ref{linear}) through the points of this curve.
\end{theo}

\subsection{Potential density}

In this subsection, we exclude the case when $D$ is a divisor.

Let $\mu:\tilde{D}\to X$ denote a desingularization; $\tilde{D}$ is
equipped with a rational self-map $\tilde{f}$ satisfying $\mu\tilde{f}=
f\mu$.

Our proof is a case-by-case analysis on the Kodaira dimension of $D$ (meaning
the Kodaira dimension of any desingularization). To start with, this 
Kodaira dimension cannot be maximal because of the presence of a self-map
of infinite order (that $f$ cannot be of finite order on $D$ is
immediate, for instance, from $f^*\sigma=-2\sigma$).  
In \cite{AV} , we already have simple geometric arguments
ruling out the cases of $\kappa(D)=-\infty$ and $\kappa(D)=0$.
The case $\kappa(D)=-\infty$ is especially simple since then the holomorphic 
2-form
would be coming from the rational quotient of $D$, but $\Sigma$ obviously
must dominate the rational quotient and this cannot be isotropic. The
case $\kappa(D)=0$ is less easy and uses the fact that $Pic(X)=\Z$ or,
equivalently, that the Hodge structure $H^2(X)^{prim}$ is irreducible
of rank 22. Namely, an argument using Minimal Model theory and the
existence of an holomorphic 2-form on $D$ gives that
$D$ must be rationally dominated by an abelian threefold or by a product of
a K3 surface with an elliptic curve. But the second
transcendental Betti number of those varieties cannot exceed 21, which
contradicts the fact that $\tilde{D}$ carries an
irreducible Hodge substructure of rank 22; see \cite{AV} for details.

Let us deal with the case $\kappa(D)=2$. We need the following lemma:

\begin{lemm}\label{orderthree} On a general $X$, the points of period 3
with respect to $f$ form a curve.
\end{lemm} 

{\it Proof:} Let $l_1$ be (a line corresponding to) such a point, $l_2=f(l_1)$,
$l_3=f^2(l_1)$, so that $f(l_3)=l_1$. There are thus planes $P_1, P_2, P_3$,
such that $P_1$ is tangent to $V$ along $l_2$ and contains $l_3$, etc.
Clearly, $P_1\neq P_2\neq P_3$. The span of the planes $P_j$ is a projective
3-space $Q$. Let us denote the two-dimensional cubic, intersection of $V$ 
and $Q$, by $W$. We can choose the coordinates $(x:y:z:t)$ on $Q$ such that
$l_1$ is given by $y=z=0$, etc. Then the intersection of $W$ and $P_1$ is
given by the equation $z^2y=0$, etc. The only other monomial from the
equation of $W$, up to a constant, can be $xyz$, since it has to be 
divisible by the three
coordinates. Therefore $W$ is a cone (with vertex at 0) over the 
cubic given by the
equation $$ax^2y+by^2z+cy^2z+dxyz=0$$ in the plane at infinity.
Now a standard dimension count (\cite{A}) shows that a general cubic admits a
one-parameter family of two-dimensional linear sections which are cones. 
Each cone on $V$ gives rise to a plane cubic on $X$. This cubic is invariant
under $f$, and $f$ acts by multiplication by $-2$ (for a suitable choice
of zero point). The points of period 3
with respect to $f$ lie on such cubic and are their points of $9$-torsion.

\begin{rema}\label{orderthree-indet} In fact the lemma says slightly more:
it applies to the indeterminacy points which are "3-periodic in the 
generalized sense",
that is, points appearing if one replaces the condition 
"$f(l_1)=l_2$, $f(l_2)=l_3$, $f(l_3)=l_1$" by "$l_2\in f(l_1)$, 
$l_3\in f(l_2)$, $l_1\in f(l_3)$"; here
by $f(l_1)$ we mean the rational curve which is the image of $l_1$ by 
the correspondence which is the graph of $f$ (equivalently, $l_2\in f(l_1)$
says that for some plane $P_3$ 
tangent to $V$ along
$l_1$, the residual line in $P_3\cap V$ is $l_2$).
\end{rema}     

\medskip    

By blowing-up $\tilde{D}$, we may assume that the Iitaka fibration 
$\tilde{D}\to B$ is regular. Its general fiber is
an elliptic curve. By \cite{NZ}, the rational self-map $\tilde{f}$ descends to $B$
and induces a transformation of finite order, so the elliptic curves
are invariant by a power of $\tilde{f}$. From proposition \ref{linear}, we
obtain that they are in fact invariant by $\tilde{f}$ itself: indeed, locally
in a neighbourhood of a fixed point, the curves invariant by $f$ are
the same as the curves invariant by its power. On a general elliptic curve, 
there is
a finite (non-zero) number of points of period three, since $\tilde{f}$ acts as
multiplication by $-2$. We have two possibilities:

\smallskip

1) These are mapped to points of period three
(in the "generalized sense" as in the \ref{orderthree-indet}) on $X$
(or the surface they form is contracted to any other curve on $X$).
Then any preimage of our surface by an iteration of $\tilde{f}$ is contracted 
as 
well, but since
there are infinitely many of them, this is impossible.

Indeed, this would be completely obvious from $\mu \circ \tilde{f}=f\circ \mu$ 
(recall that
$\mu$ denotes the map from $\tilde{D}$ to $X$) if $f$ were regular. Since $f$ is only
rational, we have to consider specially the case when either the surface $F_3$ of
period-three points on $\tilde{D}$ or any of its preimages is mapped into the
indeterminacy locus $S$. But since the map $f$ is defined along the image of 
a general elliptic
curve of $\tilde{D}$ (or at least each of its branches), 
one sees that also in this case $F_3$
must be contracted to a set of points described by \ref{orderthree-indet}; 
moreover, we know from \cite{A} that the resolution of indeterminacy of $f$
obtained by blowing-up $S$ does 
not contract surfaces, and so the preimages of $F$ must be contracted by $\mu$ as well
even if some of them are mapped to $S$.

2) This surface dominates a component of the surface of fixed points of $f$.
 In this case, several points of
period three must collapse to the same fixed point $p$. But then the
resulting branches of each elliptic curve near the generic fixed point
are interchanged by $f$, which contradicts the local description
of $f$ in \ref{linear}.

This rules out the possibility $\kappa(D)=2$.

Finally, let us consider the case $\kappa(D)=1$. The Iitaka fibration
$\tilde{D}\to C$  maps $\tilde{D}$ to a curve $C$ and the general fiber $U$ is of Kodaira
dimension 0. As before, by \cite{NZ} $\tilde{f}$ induces a finite order
automorphism on $C$, and one deduces from \ref{linear} that this
is in fact the identity. We have two possible cases:

\medskip

{\it Case 1: $U$ is not isotropic with respect to the holomorphic 2-form $\sigma$.} We use the idea from \cite{AV} as in the case $\kappa(D)=0$. 
Namely, since $X$ is generic, the Hodge structure $H^2_{prim}(X, \Q)$ is 
simple. As the 
pull-back of $\sigma$ to $\tilde{D}$ is non-zero, $H^2(\tilde{D},\Q)$ carries
a simple Hodge substructure of rang 22. Since $U$ is non-isotropic, the
same is true for $U$, but a surface of Kodaira dimension zero never
satisfies this property.

\medskip

{\it Case 2: $U$ is isotropic with respect to $\sigma$.} The kernel of the
pull-back $\sigma_D$ of $\sigma$ to $\tilde{D}$ gives a locally free subsheaf 
of rank one
in the tangent bundle $T_{\tilde{D}}$, which is in fact a subsheaf of $T_U$ 
since
$U$ is isotropic. There is thus a foliation in curves on $U$, and this
foliation has infinitely many algebraic leaves (these are intersections
of $U$ with the iterates of our original surface $\Sigma$). By Jouanolou's
theorem, this is a fibration. In other words, $D$ is (rationally) fibered 
over a 
surface $T$ in integral curves of the kernel of $\sigma_D$, and $U$ project
to curves. 
These cannot be rational curves since the surface $T$ is not uniruled
(indeed, the form $\sigma_D$ must be a lift of a holomorphic 2-form on $T$).
Therefore these are elliptic curves, and since $\kappa(U)=0$, so are
the fibers of $\pi: D\dasharrow T$.

Recall from \ref{linear} that either $D$ contains $F$, or it contains a
curve on $F$; and in this last case, locally near generic such point,
$D$ is a fibration in (isotropic) two-dimensional disks over a curve; 
in particular,
such a point is a smooth point of $D$.  If $D$ contains $F$, we get
a contradiction with \ref{fixed}: indeed, $F$ must be dominated by a union
of fibers of $\pi$, but $F$ is of general type and the fibers are elliptic. 
If $D$ contains a curve on $F$, then we look at the
"leaf" (image of $\C^2$ from the Proposition \ref{linear}) at a general point $q$ of this curve. 
If the closure of such a leaf does not coincide with the image of $U$, their 
intersection is an invariant curve, that is, the image of a line
 through the origin.
Since $U$ and the leaf are both isotropic, this must be an integral curve 
of the kernel of the restriction of $\sigma$ to $D$. But $U$
varies in a family, and this implies that the restriction of $\sigma$ to 
$D$ is zero at $q$,
a contradiction since $\sigma$ is non-degenerate.

Finally, if the image of $U$ and the closure of a general leaf coincide, then $U$ is mapped
onto itself by $f|_U$ and at some point of $U$ the tangent map is just the 
multiplication by $-2$. We recall that by definition $\kappa(U)=0$ and so
either the geometric genus of $U$ or its $m$-th plurigenus for some $m>0$ is
equal to $1$. The rational self-map $f$ preserves the spaces of pluricanonical
forms, so it must multiply some non-zero pluricanonical form by a scalar.
But because of what we know about the tangent map, this scalar can only be
equal to 4, and so the degree of $f|_U$ is 16, contradicting the calculations
of \cite{A} as in \cite{AV} (proof of Proposition 2.3).

\medskip

We thus come to a conclusion that $D$ cannot be a divisor, so $D=X$.

\begin{rema}\label{useofgenericity} Here, unlike in the proof of
non-preperiodicity, we do assume that $X$ is "sufficiently general"
and so $Pic(X)=\Z$. It would be interesting to check whether one
can modify the argument to get rid of this assumption.
\end{rema}

\section{Appendix: a version of Siegel's theorem}

In this appendix we explain how to modify the proof of Siegel's 
theorem on linearization of $p$-adic diffeomorphisms given in 
\cite[Theorem 1, \S4, p.423]{HY} in order to adapt it to the 
situation where the fixed point is not isolated.

Let $k$ be a complete non-archimedian field and $n>r\ge 0$ be integers. 

\begin{theo} \label{linearizaciq}
Let $f=(f^{(1)},\dots,f^{(n)})$ be an analytic diffeomorphism of an 
$n$-dimensional 
domain over $k$. Assume that the fixed set of $f$ is $r$-dimensional 
and that the tangent maps of $f$ are semisimple at all fixed points of $f$. 
Suppose, moreover, that the eigenvalues of $Df$ distinct from $1$ (at fixed 
points of $f$) are either 
\begin{enumerate} \item equal and are not roots 
of unity at a general fixed point, or \item constant and, if denoted by 
$\lambda_{r+1},\dots,\lambda_n$, satisfy the following 
bad diophantine approximation property: 
$|\lambda_{r+1}^{i_{r+1}}\cdots\lambda_n^{i_n}-\lambda_j|_k\ge 
C(i_{r+1}+\dots+i_n)^{-\beta}$ for some $C>0$ and $\beta\ge 0$, 
any $r<j\le n$ and $(i_{r+1},\dots,i_n)\in\Z^{n-r}_{\ge 0}$ 
such that $i_{r+1}+\dots+i_n\ge 2$. \end{enumerate}
Then there exist coordinates $x_1,\dots,x_n$ in a neighbourhood of a 
general fixed point of $f$ such that 
$f^{(i)}(x)=\lambda_i(x_1,\dots,x_r)x_i$ for all $1\le i\le n$, 
where $\lambda_1=\dots=\lambda_r=1$ (and in the second case 
$\lambda_i(x_1,\dots,x_r)$ are constant for all $1\le i\le n$). \end{theo}

\vspace{4mm} 

The desired coordinates are constructed, following the Newton's method, 
as the limit (in an appropriate non-archimedian Banach space) of a 
certain sequence 
of approximations $h_0,h_1,h_2,\dots$, which are diffeomorphisms of the open 
$\rho_{\infty}$-neighbourhood of a general fixed point of $f$ for some 
$\rho_{\infty}>0$. (In fact, $h_i$ is also a diffeomorphism of the 
open $\rho_i$-neighbourhood for each $i\ge 1$, where 
$\rho_0>\rho_1>\rho_2>\dots>\rho_{\infty}$ is a decreasing sequence of radii.) 

After a preliminary step (Corollaries \ref{equal-eigenv}, \ref{const-eigenv}), 
the procedure is the same as in \cite{HY}, but instead of working in the 
spaces $A^2_{\rho}(k^n)$ and $B^2_{\rho}(k^n)$ of loc. cit. we work in 
smaller (when $r>0$) spaces $A^{(r)}_{\rho}(k^n)$ and $B^{(r)}_{\rho}(k^n)$, 
cf. below.

\begin{lemm}\label{blok-diag} 
Let $f$ be an analytic diffeomorphism of an $n$-dimensional domain over 
$k$, $F$ be the fixed set of $f$, and $q$ be a general point of $F$. 
Assume that $F$ is $r$-dimensional and that only $r$ eigenvalues of 
the tangent map of $f$ at $q$ are equal to $1$. Then there exist 
coordinates $x_1,\dots,x_n$ in a neighbourhood of $q$ such that 
$x_{r+1},\dots,x_n$ generate the ($f$-invariant) ideal $I_F$ of 
$F$ and $f^{(i)}(x)\equiv x_i\pmod{I_F^2}$ for all $1\le i\le r$. 
\end{lemm}
{\it Proof.} For a multi-index $I=(i_1,\dots, i_k)$, denote 
$|I|=\sum_{j=1}^{j=k}i_j$.

We split the coordinates into two groups: 
$x':=(x_1,\dots,x_r)$ and $x'':=(x_{r+1},\dots,x_n)$. Then 
$f(x)=x+\sum_{I\in\Z_{\ge 0}^{n-r}:~|I|\ge 1}a_I(x')(x'')^I$ for some 
analytic $a_I(x')$, so $f(x)\equiv(x'+\sum_{I\in\{0,1\}^{n-r}:~|I|=1}
a_I'(x')(x'')^I;x''+\sum_{I\in\{0,1\}^{n-r}:~|I|=1}a_I''(x')(x'')^I)
\pmod{I_F^2}$, so $f\equiv(x'+a'x'';x''+a''x'')\pmod{I_F^2},$ 
where $a'=a'(x')$ is the matrix with columns $a'_I(x')$ for $|I|=1$ 
and $a''=a''(x')$ is the matrix with columns $a''_I(x')$ for $|I|=1$
(so $a''$ is an $(n-r)\times (n-r)$-matrix, invertible since $Df$ is
invertible at $q$). 
Let $h(x)=(x'+a'(a'')^{-1}x'';x'')$. Then 
$h^{-1}(x)\equiv(x'-a'(a'')^{-1}x'';x'')\pmod{I_F^2}$, 
$f(h(x))\equiv(x'+a'(a'')^{-1}x''+a'x'';x''+a''x'')=
(x'+a'(1+(a'')^{-1})x'';x''+a''x'')$, and finally,
$h^{-1}(f(h(x)))\equiv(x';x''+a''x'')\pmod{I_F^2}$. 

\begin{coro} \label{equal-eigenv} In the setting of Lemma \ref{blok-diag}, 
assume that $({\rm i})$ 
the tangent maps of $f$ are semisimple at all points of $F$ and $({\rm ii})$ 
their eigenvalues distinct from $1$ are equal. Then there exist coordinates 
$x_1,\dots,x_n$ in a neighbourhood of a general point of $F$ such that 
$x_{r+1},\dots,x_n$ generate the ideal $I_F$ of $F$, 
$f^{(i)}(x)\equiv x_i\pmod{I_F^2}$ for all $1\le i\le r$ and 
$f^{(i)}(x)\equiv\lambda(x')x_i\pmod{I_F^2}$ for all $r<i\le n$. \end{coro}
{\it Proof.} In the setting of the proof of Lemma \ref{blok-diag}, 
the matrix $a''$ is diagonalizable with the same eigenvalues, 
so it is already diagonal. 

\begin{coro}\label{const-eigenv}  In the setting of Lemma \ref{blok-diag}, 
assume that $({\rm i})$ 
the tangent maps of $f$ are semisimple at all points of $F$ and $({\rm ii})$ 
their eigenvalues $\lambda_1,\dots,\lambda_n$ do not vary. 

Then there exist coordinates $x_1,\dots,x_n$ 
in a neighbourhood of a general point of $F$ such that 
$x_{r+1},\dots,x_n$ generate the ideal $I_F$ of $F$ and 
$f^{(i)}(x)\equiv\lambda_ix_i\pmod{I_F^2}$ for all $1\le i\le n$. \end{coro}
{\it Proof.} In the setting of the proof of Lemma \ref{blok-diag}, there 
is transform of the coordinates $x_1,\dots,x_n$, identical on $x_1,\dots,x_r$, 
which is linear on $x_{r+1},\dots,x_n$ with coefficients in functions of 
$x_1,\dots,x_r$ making the matrix diagonal. 

Namely, after a $k$-linear change of variables we can assume that 
$a''$ is diagonal at $q$. Let 
$p_i:=\prod_{j:~\lambda_j\neq\lambda_i}(\lambda_i-\lambda_j)^{-1}
(a''-\lambda_j)$ be a projector onto the $\lambda_i$-eigenspace of $a''$. 
If $\{e_{r+1},\dots,e_n\}$ is an eigenbasis of $a''(q)$, considered as 
sections of the restriction of the tangent bundle to $F$, then 
$\{E_i:=p_ie_i\}_{r<i\le n}$ is a system of eigenvectors of $a''$ and its 
reduction modulo the maximal ideal in $k[[x_1,\dots,x_r]]$ is the eigenbasis 
$\{e_{r+1},\dots,e_n\}$ of $a''(q)$. This means that $E_1,\dots,E_n$ 
generate the tangent bundle, cf. \cite[Proposition 2.8]{AM}. Then the dual 
basis of the cotangent bundle gives the desired linear transformation of 
the coordinates $x_{r+1},\dots,x_n$ with coefficients in functions of 
$x_1,\dots,x_r$. 

\vspace{4mm}

Now we modify \cite{HY} to the setting of Corollary \ref{const-eigenv} 
(to allow some ``resonances''). 
Let $f$ and $x_1,\dots,x_n$ be as in Corollary \ref{const-eigenv} and 
$\Lambda$ be the differential of $F$ at $q$, i.e., an $n\times n$ diagonal 
matrix with entries $\lambda_1,\dots,\lambda_n$, the first $r$ ones we 
assume to be equal to 1. We assume that $\lambda_{r+1},\dots,\lambda_n$ 
satisfy the following bad diophantine approximation property: 
$|\lambda_{r+1}^{i_{r+1}}\cdots\lambda_n^{i_n}-\lambda_j|_k\ge 
C(i_{r+1}+\dots+i_n)^{-\beta}$ for some $C>0$ and $\beta\ge 0$, 
any $r<j\le n$ and $i=(i_{r+1},\dots,i_n)\in\Z^{n-r}$ such that 
$i_{r+1}+\dots+i_n>1$. 

For real $\rho>0$ define the spaces 
\begin{itemize} 
\item \label{Banach-algebra} $A_{\rho}(k^n):=\{\phi=\sum_Ia_Ix^I\in 
k[[x_1,\dots,x_n]]~|~\sup|a_I|\rho^{|I|}=:\|\phi\|_{\rho}<\infty\}$ 
(this is a non-archimedian Banach algebra);
\item $A^{(r)}_{\rho}(k^n)$ as the ideal in $A_{\rho}(k^n)$ generated by 
$x_ix_j$ for all $r<i,j\le n$ 
($A^{(0)}_{\rho}(k^n)$ is denoted in \cite{HY} by $A^2_{\rho}(k^n)$); 
\item $B^{(r)}_{\rho}(k^n):=\{\phi\in A^{(r)}_{\rho}(k^n)
~|~\phi\circ\Lambda\in A^{(r)}_{\rho}(k^n)\}$. In particular, 
$B^{(r)}_{\rho}(k^n)=A^{(r)}_{\rho}(k^n)$ if 
$|\lambda_1|=\dots=|\lambda_n|=1$. \end{itemize}

The set $x+(A^{(r)}_{\rho}(k^n))^n$ is a group (with respect to the 
composition). This group acts on $A^{(r)}_{\rho}(k^n)$. We consider 
$f$ as an element of $A_{\rho}(k^n)^n$ for an appropriate $\rho$. In fact,
since $f$ satisfies the assumptions of Corollary 1.4, 
$f-\Lambda\in (A^{(r)}_{\rho}(k^n))^n$.

We are looking for a root of the equation $F_f(h):=f\circ h-h\circ\Lambda=0$ 
in the form $h(x)=\lim\limits_{i\to\infty}h_i(x)$, where 
$h_i(x)-x\in(A^{(r)}_{\rho_i}(k^n))^n\subset(A^{(r)}_{\rho_{\infty}}(k^n))^n$ 
for some $\rho_0>\rho_1>\dots>\rho_{\infty}>0$. Replacing if necessary $f$ 
by $uf(u^{-1}x)$ ($u\in k$, $|u|_k\gg 1$), we may assume that 
$f\in(A_1(k^n))^n$ and $\|f-\Lambda\|_1$ is as small as we want. 
Let $\rho_i=1/2+2^{-i-1}$, $\rho_{\infty}=1/2$ for all integer $i\ge 0$. 

The sequence $(h_i)$ is constructed inductively, following the Newton's 
method. Let $h_0=id$ and $h_{i+1}=h_i+\Delta_i$. Here we have to define
$\Delta_i$ in order to make the induction step. One has 
\begin{multline*}F_f(h_i+\Delta_i)=f\circ(h_i+\Delta_i)-
f\circ h_i+F_f(h_i)-\Delta_i\circ\Lambda\\ =[f\circ(h_i+\Delta_i)-f\circ h_i+
DF_f(h_i)\cdot E_i-Df\circ h_i\cdot\Delta_i]+F_f(h_i)-(Dh_i\circ\Lambda)LE_i,
\end{multline*} 
where $E_i:=(Dh_i)^{-1}\cdot\Delta_i$ and 
$L:(B^{(r)}_{\rho}(k^n))^n\to(A^{(r)}_{\rho}(k^n))^n$ is the linear operator, 
defined by $Lw=w\circ\Lambda-\Lambda\circ w$. 

We define $E_i$ (or equivalently, $\Delta_i$) so that 
$F_f(h_i)=(Dh_i\circ\Lambda)LE_i$. To do this, we have to invert $L$.
The injectivity of $L$ is evident: $\phi=\sum_Ia_Ix^I\mapsto(\sum_I(\lambda^I-
\lambda_i)a^{(i)}_Ix^I)_{1\le i\le n}$, where $\lambda^I=\lambda_1^{i_1}
\cdots\lambda_n^{i_n}$. Moreover, it is also evident that for any 
$g\in(A^{(r)}_{\rho}(k^n))^n$ there exists a unique $n$-tuple of formal 
series $w$ such that $Lw=g$, $w(x_1,\dots,x_r,0,\dots,0)=0$ and $\frac
{\partial w}{\partial x_i}(x_1,\dots,x_r,0,\dots,0)=0$ for all $r<i\le n$. 

It remains to show that $\|\Delta_i\|_{\rho_{\infty}}$ tends to 0. 
But this is explained in \cite[\S4.4]{HY}, if one replaces \cite[Lemma 15]{HY} 
by the following version (with the same proof):

\begin{lemm} \label{lemma-15-HY} For any $\delta>0$ 
and any $g\in(A^{(r)}_{\rho}(k^n))^n$ the solution $w$ of the linear 
equation $Lw=g$ belongs to $(B^{(r)}_{\rho-\delta}(k^n))^n$ and satisfies 
$\|w\|_{\rho-\delta}\le C_1\frac{\|g\|_\rho}{\delta^{\beta}}\rho^{\beta}$, 
$\|Dw\|_{\rho-\delta}\le C_1\frac{\|g\|_{\rho}}{\delta^{\beta}}
\frac{\rho^{\beta}}{\rho-\delta}$, $\|Dw\circ\Lambda\|_{\rho-\delta}\le 
C_1\frac{\|g\|_{\rho}}{\delta^{\beta}}\frac{\rho^{\beta}}{\rho-\delta}$, 
where $C_1$ is a constant depending only on $C$, $\beta$, $\|\Lambda\|$. 
{\rm (We put on $(B^{(r)}_{\rho}(k^n))^n$ the max norm: $\|\phi\|_{\rho}
=\max(\|\phi\|_{\rho},\|\phi\circ\Lambda\|_{\rho})$.)} \end{lemm}

E. Amerik: Universit\'e Paris-Sud, Laboratoire de Math\'ematiques, Campus
Scientifique d'Orsay, B\^atiment 425, 91405 Orsay, France. 

Ekaterina.Amerik@math.u-psud.fr.

\medskip

F. Bogomolov: Courant Institute of Mathematical Sciences, 251 Mercer Str.,
New York NY 10012 USA.

bogomolo@cims.nyu.edu

\medskip

M. Rovinsky: Independent University of Moscow, 119002 Moscow,
B.Vlasievsky Per. 11, and
Institute for Information Transmission Problems of Russian Academy of Sciences.

marat@mccme.ru


\begin{thebibliography}{}

\bibitem[A]{A} E. Amerik, A computation of invariants of a rational 
self-map, Ann. Fac. Sci. Toulouse 18 (2009), to appear.
\bibitem[AC]{AC} E. Amerik, F. Campana, Fibrations m\'eromorphes sur certaines vari\'et\'es \`a fibr\'e canonique trivial,  Pure Appl. Math. Q.  4  (2008),  no. 2 (special issue in honour of F. Bogomolov), part 1, 509--545.
\bibitem[AV]{AV} E. Amerik, C. Voisin, Potential density of rational 
points on the variety of lines of a cubic fourfold,  Duke Math. J. 
\textbf{145}  (2008),  no. 2, 379--408
\bibitem[Arn]{Arn} V. I. Arnold, {\bf Geometrical methods in the theory 
of ordinary differential equations.} Second edition, Grundlehren der 
Mathematischen Wissenschaften, 250, Springer-Verlag, New York, 1988. 

 \bibitem[AtM]{AM} M.F.Atiyah, I.G.Macdonald, {\bf Introduction to commutative 
algebra.} Addison-Wesley Publishing Co., Reading, Mass.-London-Don Mills, 
Ont. 1969. 
\bibitem[Bv]{Bv} A. Beauville, Some remarks on K\"ahler manifolds with $c_{1}=0$, in: Classification of algebraic and analytic manifolds (Katata, 1982),  
1--26, Progr. Math., 39, Birkhäuser Boston, Boston, MA, 1983.
\bibitem[BT]{BT} F. Bogomolov, Yu. Tschinkel, Density of rational points
 on elliptic $K3$ surfaces,  Asian J. Math.  4  (2000),  no. 2, 351--368.
\bibitem[C]{C} F. Campana, Orbifolds, special varieties and classification theory,  Ann. Inst. Fourier (Grenoble)  54  (2004),  no. 3, 499--630. 
\bibitem[CG]{CG} H. Clemens, Ph. Griffiths, The intermediate Jacobian of the 
cubic threefold,  Ann. of Math. (2)  95  (1972), 281--356. 

\bibitem[E]{eisen} D.Eisenbud, {\bf Commutative algebra. With a view toward 
algebraic geometry.} Graduate Texts in Mathematics, 150. Springer-Verlag, 
New York, 1995. 
\bibitem[F]{F} N. Fakhruddin, Questions on self-maps of algebraic varieties,
J. Ramanujan Math. Soc. 18 (2003), no. 2, 109--122. 
\bibitem[GT]{GT} D. Ghioca, T. Tucker: Periodic points, linearizing maps, 
and the dynamical Mordell-Lang problem,  J. Number Theory  129  (2009),  
no. 6, 1392--1403.
\bibitem[HT1]{HT} B. Hassett, Yu. Tschinkel, Abelian fibrations and 
rational points on symmetric products,
Internat. J. Math. 11 (2000), no. 9, 1163--1176.
\bibitem[HT2]{HT2} B. Hassett, Yu. Tschinkel, Flops on holomorphic symplectic fourfolds and determinantal cubic hypersurfaces, preprint arXiv:0805.4162. 
\bibitem[HY]{HY} M. Herman, J.-C. Yoccoz, Generalizations of some theorems of
small divisors to non-Archimedean fields, in Geometric dynamics 
(Rio de Janeiro, 1981),  408--447, Lecture Notes in Math., 1007, Springer, 
Berlin, 1983. 


\bibitem[NZ]{NZ} N. Nakayama, D.-Q. Zhang: Building blocks of \'{e}tale 
endomorphisms of complex projective manifolds, RIMS preprint no. 1577, 2007.
\bibitem[Pc]{Pc} G. Pacienza, Rational curves on general projective 
hypersurfaces,  J. Algebraic Geom.  12  (2003),  no. 2, 245--267.
\bibitem[P]{P} H. Poincar\'e, {\bf Oeuvres.} T. I, p. XXXVI-CXXIX, Gauthier-Villars, Paris, 1928.
\bibitem[Rg]{Rg} F. Rong, Linearization of holomorphic germs with 
quasi-parabolic fixed points, Ergodic Theory Dynam. Systems  28  (2008),  
no. 3, 979--986. 
\bibitem[T]{T} T. Terasoma, Complete intersections with middle Picard 
number 1  defined over ${\mathbb Q}$, Math. Z. 189 (1985), 289--296.

\bibitem[V]{V} C. Voisin, A correction: ``On a conjecture of Clemens on rational curves on hypersurfaces'', J. Differential Geom.  49  (1998),  no. 3, 601--611. 
\bibitem[V1]{V1} C. Voisin, Intrinsic pseudo-volume forms and $K$-correspondences, in: The Fano Conference,  761--792, Univ. Torino, Turin, 2004.


\bibitem[Yu]{Yu} K. Yu, Linear forms in $p$-adic logarithms. II, 
Compositio Math. 74 (1990), no. 1, 15--113. 

\bibitem[Z]{Z} S.-W. Zhang, Distributions in algebraic dynamics, in: Surveys 
in differential geometry, Vol. X,  381--430, Int. Press, Somerville, MA, 2006. 


\end{thebibliography}
\end{document}